\documentclass[]{article}

\usepackage{amssymb}
\usepackage{ART-E2VEM-Mixed}
\usepackage{ART-E2VEM-Mixed-Bis}
\usepackage[margin=1.5in]{geometry}
\usepackage{comment}

\title{A lowest order stabilization-free mixed Virtual Element Method}
\author{Andrea Borio, Carlo Lovadina, Francesca Marcon and Michele Visinoni}
\date{}

\begin{document}
\maketitle

\begin{abstract}
We initiate the design and the analysis of stabilization-free Virtual Element Methods for the laplacian problem written in mixed form. A Virtual Element version of the lowest order Raviart-Thomas Finite Element is considered. To reduce the computational costs, a suitable projection on the gradients of harmonic polynomials is employed. A complete theoretical analysis of stability and convergence is developed in the case of quadrilateral meshes. Some numerical tests highlighting the actual behaviour of the scheme are also provided.  
\end{abstract}

\section{Introduction}
In these years, the study of numerical methods for solving partial differential equations on polygonal/polytopal meshes has been experiencing a growing interest in the scientific community.
In particular, one of the most recent developments in this field is represented by the Virtual Element Method (VEM).
This technology was first introduced in the primal conforming Poisson problem in~\cite{Beirao2013a} as a generalization of $\sobh{1}{}$-conforming Finite Element Method. Successively, the extension to the $\mathrm{H}(\div)$-conforming vector fields, generalizing Mixed Finite Elements~\cite{BrezziBoffiFortin}, has been introduced in \cite{Brezzi2014} and developed in \cite{Beirao2016,Beirao2016b,Dassi2020}.
Thanks to the great flexibility of the method, both primal and mixed formulation of VEM have been applied to a large range of applications, such as elastic and inelastic problems \cite{Beirao2015a,Artioli2017,Dassi2020b,Dassi2021b}, simulations in fractured media \cite{Benedetto2016b,Benedetto2016c,BENEDETTO2022114204,Berrone2023} and in porous media mechanics \cite{Borio2011,Berrone2021,Berrone2022b}, just to mention a few of them.

The key ideas of VEM may be summarised as follows.

\begin{itemize}
\item The local approximation spaces are defined as the solutions to suitable local partial differential problems. Therefore, VEM functions are not explicitly known, but only a limited information is available. However, the local approximation spaces contain polynomials up to a suitable degree. 
\item A computable projection onto a polynomial space is involved. Typically, the projection is valued onto the polynomials contained in the approximation spaces. 
\item The discrete bilinear forms are characterized by the sum of a singular part maintaining consistency on polynomials, and a stabilizing form enforcing coercivity. 
\end{itemize}

However, in general the stabilising form mentioned above is designed without a clear physical meaning, but only requiring minimal assumptions to make the method stable. Though efficient recipes to tune the stabilisation term have been proposed (see for instance \cite{d-recipe1,d-recipe2}), in certain complex situations it might be preferable to avoid dealing with the choice of such forms. As examples, we mention highly non-linear problems; problems where highly anisotropic meshes occurs; advection-diffusion problems. In addition, the stabilization term could be problematic in connection with the analysis of a-posteriori error estimates \cite{Cangiani2017,Berrone2017} (however, the recent work \cite{Canuto2022StabFreeArxiv} presents a first study which provides stabilization-free upper and lower a-posteriori bounds for triangular meshes with hanging nodes).

Virtual Element schemes for which no stabilisation form is required have been recently presented, in different 2D frameworks, in \cite{BBME2VEM,Berrone2022,DAltri2021,Lamperti2022}. These approaches share the idea to employ a projection onto a polynomial space of higher degree than the one usually taken in standard VEM.
It is worth noticing that the polynomial degree depends on the number of edges of each polygon: as expected, it increases as the edge number gets larger. As a consequence, the quadrature computational cost significantly grows in presence of elements with many edges, without any improvement in the convergence rate.

This paper follows similar lines of the above-mentioned stabilisation-free attempts \cite{BBME2VEM,Berrone2022,DAltri2021,Lamperti2022}, but for the Laplacian problem written in the usual $H(\div)-L^2$ mixed formulation.
In particular, we consider a VEM version of the lowest order Raviart-Thomas Finite Element Method, see \cite{Beirao2016}. To reduce the computational cost connected to quadrature, a suitable projection operator onto the gradients of \emph{harmonic} polynomials is selected, similarly to the scheme introduced for the primal formulation in \cite{BBMTLetter}. The resulting scheme has the following features.

\begin{itemize}
\item It is a conforming mixed VEM method for which no stabilization term is needed. 
\item The method shows first order convergence rate for the natural norms and, in most cases, a behaviour comparable with the standard lowest order Raviart-Thomas VEM for which the stabilisation term is suitably tuned. However, for highly anisotropic meshes, our method seems to display a better performance.     
\item Despite a projection over higher-order polynomial spaces is employed, the use of \emph{harmonic} polynomials greatly alleviate the additional computational costs. 
\end{itemize}
These properties indicate that the present approach could be a valid alternative to the lowest order Raviart-Thomas Virtual Element Methods, especially in those complex situations where, for the latter scheme, a particular care in the treatment of the stabilising form is required.

From a theoretical point of view, the present paper can be considered as a first contribution, since we present a rigorous analysis only for the quadrilateral case (of course, the similar arguments could be applied also for triangular elements). However, the general theory for polygons with an arbitrary number of edges is not currently available and will be treated in a future work.

A brief outline of the paper is as follows. In Section \ref{sec:model} we define the model problem. 
Section \ref{sec:discr_formulation} contains the statement of the discrete problems, introducing all the bilinear and linear forms involved.
In section \ref{sec:well-posedness}, we prove the well-posedness of the discrete problem in the quadrilateral case. For the same kind of meshes, we derive optimal error estimates in Section \ref{sec:apriori} and, finally, in Section \ref{sec:numtests} we present some numerical results that assess the convergence rate of the method; a comparison with the standard lowest order Raviart-Thomas VEM is also provided.

\section{Model problem}\label{sec:model}
Let $\Omega\subset\mathbb{R}^2$ be a computational domain. We are interested in studying the following mixed formulation of the Poisson problem:
\begin{equation}\label{eq:poissonMixed}
\begin{cases}
-\div\bs\sigma=f \quad &\text{in} \; \Omega \\
\bs\sigma = \nabla u\quad &\text{in} \; \Omega \\
u = 0 \quad &\text{on} \; \partial\Omega
\end{cases},
\end{equation}
where the forcing term $f\in\lebl{\Omega}$. We consider homogeneous natural boundary conditions only for sake of simplicity: the extension to non-homogeneous or essential boundary conditions can be treated with the same techniques used for other more classical Galerkin methods, such as the FEM.
Let $\scal[\Omega]{\cdot}{\cdot}$ denote the $\lebl{}$ scalar product and $a(\bfsigma,\bftau):=\scal[\Omega]{\bfsigma}{\bftau}$, then the mixed variational formulation of \eqref{eq:poissonMixed} is given by: find $\left(\bs \sigma, u\right)\in \Sigma\times U$, where $\Sigma:=\sobhDiv{\Omega}$ and $U:=\lebl{\Omega}$ such that
\begin{equation}\label{eq:cont-var-form}
\begin{cases}
a(\bfsigma,\bftau)+\scal[\Omega]{\div\bs\tau}{u} = 0 \quad &\forall \; \bs\tau\in\Sigma\,,\\
\scal[\Omega]{\div\bs\sigma}{v} = -\scal[\Omega]{f}{v} \quad &\forall \; v\in U.
\end{cases}
\end{equation}
Well posedness of the above problem \eqref{eq:cont-var-form} is standard and the details can be found, for instance, in \cite{BrezziBoffiFortin}.
\section{VEM discrete formulation}
\label{sec:discr_formulation}
%
In order to state the discrete formulation of \eqref{eq:cont-var-form}, let $\Mh$ be a polygonal tessellation of $\Omega$. For every element $E\in\Mh$, its area and diameter are denoted by $|E|$ and $h_E$, respectively. 
As usual, the maximum of the diameters $h_E$ for $E\in\Mh$ is the mesh size, denoted by $h$, i.e. $h= \max_{E\in\Mh} h_E$.
We assume that each $E\in\Mh$ is such that
\begin{enumerate}[label=\textbf{A.\arabic*}]
\item\label{meshA_1} $E$ is star-shaped with respect to a ball of radius $\geq\gamma h_E$,
\item\label{meshA_2} for any edge $e$ of $\partial E$, $\abs{e}\geq \gamma h_E$,
\end{enumerate}  
where $\gamma$ is a positive constant.

To continue, for any given $E\in \Mh$ and non-negative integer $k$, $\Poly{k}{E}$ denotes the space of polynomials of degree up to $k$ defined on $E$. Moreover, we introduce $\PolyH{k}{E}\subseteq \Poly{k}{E}$ as the space of \emph{harmonic} polynomials of degree up to $k$ defined on $E$; the dimension of this latter space is $2k+1$.  

\subsection{The local spaces}\label{subsec:localSpaces}
In this section we introduce the discrete local space and their interpolation properties.
Given a generic quadrilateral $E\in\Mh$, we introduce the following local VEM space:
\begin{equation}\label{eq:vemDiscreteSpace}
\begin{split}
\Vh[E] &:= \left\{\bs{\tau}_h\in\sobhDiv{E}\colon \exists v \in
\sobh{1}{E}\,\mbox{ s.t. }\,\bs{\tau}_h = \nabla v\,,\right.
\\
&\left. \,\qquad \bs\tau_h\cdot\bbn_e \in \Poly{0}{e} \quad \forall e
\in{\partial E}\,, \quad \div\bs{\tau}_h\in\Poly{0}{E} \right\}.
\end{split}
\end{equation}
Accordingly, for the local space $\Vh[E]$ the following degrees of freedom can be taken:
\begin{equation}\label{eq:defDofSigma}
	\bs\tau_h \rightarrow \frac{1}{|e|}\int_{e} \bs\tau_h\cdot\bbn_{e}~\de = \bs\tau_h\cdot\bbn_{e}, \qquad \forall e\in\partial E.
\end{equation}
The unisolvence of the above degrees of freedom is proved, e.g., as in~\cite{Beirao2016b},
so that $\dim(\Vh(E)) = 4$.
We remark that, once $\bs\tau_h\cdot\bbn_e = c_e\in\Poly{0}{E}$ is given for all $e\in\partial E$, the quantity $\div \bs\tau_h\in\Poly{0}{E}$ is uniquely determined. 
Since $\div\bs\tau_h\in\Poly{0}{E}$ then
\begin{equation}
	\div\bs\tau_h = \frac{1}{|E|}\int_E\div\bs\tau_h~\dEl =\frac{1}{|E|}\sum_{e\in\partial E}\int_e \bftau_h\cdot\bbn_e~\de= \frac{1}{|E|}\sum_{e\in\partial E} |e|c_e.
\end{equation}
The local approximation space for $U$ is simply defined as follows
\begin{equation}~\label{eq:polynomialDiscreteSpace}
\Uh(E):=\left\{u_h\in L^2(E): u_{h}\in \Poly{0}{E}\right\}.
\end{equation}
Accordingly, for the local space $\Uh(E)$ the following degrees of freedom can be taken:
\begin{equation}
	u_h \rightarrow \frac{1}{|E|}\int_E u_h~\dEl.
\end{equation}
It follows that $\dim(\Uh(E)) = 1$.

\subsection{Approximation in $\Vh$ and $U_h$}\label{subsec:interpolation}
Let us consider the space $W(\Omega)= \sobhDiv{\Omega}\cap\left[\lebl[r]{\Omega}\right]^2$ ($r>2$), equipped with the natural norm.   
We define an interpolation operator 
\begin{equation}
\Ih : W(\Omega) \longrightarrow \Vh    
\end{equation}
by requiring 
\begin{equation}\label{eq:defInterptau}
\int_e \left(\bs{\varsigma}-\Ih\bs{\varsigma}\right)\cdot\bbn_{e}~\de =  0, \quad \forall \mbox{ edge } e \mbox{ of the elements in } \Mh\,.
\end{equation}
Using the unisolvence of the degrees of freedom, e.g. see~\cite{Beirao2016b}, it is not difficult to check that such a $\Ih\bs{\varsigma}$ exists and it is unique in $\Sigma_h$.
This definition implies that for each $E\in\Mh$ 
\begin{equation}\label{eq:divCompat}
\int_E \div \left( \bs{\varsigma}-\Ih\bs{\varsigma}\right)~\dEl = 0\,.
\end{equation}
Hence, since for each $E\in\Mh$ $\div\Ih\bs{\varsigma}\in\Poly{0}{E}$, we obtain the \emph{commuting diagram property}
\begin{equation}\label{eq:propertyProjdivInterp}
\div\Ih\bs{\varsigma} = \proj{0}{E} \div\bs{\varsigma},
\end{equation}
where $\proj{0}{E}:\lebl{E}\to\Poly{0}{E}$ is the $\lebl{}$ projection operator onto constants.
We now remark that $(\Ih\bs{\varsigma})_{|E}= \nabla\varphi^\ast$, $\varphi^\ast$ being the solution to the local (compatible) Neumann problem 

\begin{equation}\label{eq:localinterp}
\begin{cases}
\Delta\varphi^\ast= \proj{0}{E} \div\bs{\varsigma} \quad &\text{in} \; E \\
\nabla \varphi^\ast \cdot\bbn_e = \proj{0}{e}(\bs\varsigma \cdot\bbn_e) \quad &\text{on every $e$ side of} \; \partial E ,
\end{cases}
\end{equation}
where $\proj{0}{e}$ denotes the $\mathrm{L}^2$ projection operator onto the constant functions on $e$.
Regularity results of elliptic equations and Sobolev embedding theorems shows that there exists $r^\ast>2$ such that for $r\in (2,r^\ast]$ it holds
\begin{equation}\label{eq:interpCont}
\norm[0,E]{\Ih\bs{\varsigma}} \leq C_{r^\ast}  ||{\bs{\varsigma}}||_{W(E)}\, .
\end{equation}
Moreover assuming $\bs{\varsigma}\in\left[\sobh{1}{\Omega}\right]^2$ and $\div\bs{\varsigma}\in\sobh{1}{\Omega}$,  
the following approximation results hold: for each $h$, for each $E\in\Mh$
\begin{equation}\label{eq:interpEstimateDivTau}
\norm[0,E]{\div(\bs{\varsigma}-\Ih\bs{\varsigma})} \leq C_d h^s_E \seminorm[s,E]{\div\bs{\varsigma}}, \quad s=0,1
\end{equation}
and
\begin{equation}\label{eq:interpEstimateTau}
\norm[0,E]{\bs{\varsigma}-\Ih\bs{\varsigma}} \leq C_\varsigma h_E \seminorm[1,E]{\bs{\varsigma}}\, .
\end{equation}
Above, $C_{r^\ast}$, $C_d$ and $C_\varsigma$ are positive constants depending only on the constant $\gamma$ of the mesh assumptions \ref{meshA_1} and \ref{meshA_2}.

Moreover, we recall that, given $w\in\sobh{1}{\Omega}\cap \mathrm{L}^2(\Omega)$, for its $\mathrm{L}^2$ projection $\proj{0}{E}w\in\Uh$ it holds for each $h$, for each $E\in\Mh$
\begin{equation}\label{eq:interpEstimateU}
\norm[0,E]{w-\proj{0}{E}w} \leq C h_E^s \seminorm[s,E]{w}, \quad s=0,1\,,
\end{equation}
where $C>0$ depends only on the constant $\gamma$ of the mesh assumptions \ref{meshA_1} and \ref{meshA_2}.

\subsection{The local forms}\label{subsec:localForms}
In this section we introduce the VEM counterparts of the local forms associated with the continuous problem.
\paragraph{The local mixed term} Given $E\in \Mh$, we notice that the term
\begin{equation*}
\left(\div \bs\tau_h,v_h\right)_{E}=\int_E v_h\div \bs\tau_h~\dEl
\end{equation*}
is computable for every $\bftau_h\in\Sigma_h(E)$ and $v_h\in U_h(E)$ via degrees of freedom. For this reason, we do not need to introduce any approximation of the continuous terms $(\div \bftau, u)$ and $(\div \bfsigma, v)$ in problem~\eqref{eq:cont-var-form}.
\paragraph{The local bilinear form $a^{E}(\cdot,\cdot)$}
The local bilinear form
\begin{equation*}
a^E(\bfsigma_h,\bs\tau_h)  = \int_E \bfsigma_h \cdot \bftau_h~\dEl
\end{equation*}
is not computable for a general pair
$(\bfsigma_h,\bftau_h)\in \Sigma_h(E)\times \Sigma_h(E)$.  Here, instead of using
the standard VEM procedure (cf. \cite{Beirao2016}), we introduce a
local self-stabilized discrete bilinear form.  Let
\begin{equation}\label{eq:defProj0}
\projhat{k-1}{E}:\lebldouble{E}\to\nabla\PolyH{k}{E}
\end{equation}
be the $\lebl{E}$-projection operator onto the space $\nabla\PolyH{k}{E}$, i.e. the space of gradients of \emph{harmonic} polynomials of degree at most $k$, with $k\ge 1$. More precisely, $\projhat{E}{}$ is
defined by the orthogonality condition: for each  $\bs\tau\in\lebldouble{E}$, it holds
\begin{equation}~\label{eq:defProj}
\scal[E]{\projhat{k-1}{E}\bs\tau}{\nabla p} = \scal[E]{\bs\tau}{\nabla p}, \qquad \quad \forall \; p\in \PolyH{k}{E}.
\end{equation} 
In order to attain stability, the approximation of $a^E(\cdot,\cdot)$ depends on the number of edges of $E$, denoted by $n_E$. More precisely, $[\cdot]$ being the integer part, we select
\begin{equation}\label{k-choice}
k=\left[\frac{n_E+1}{2}\right]    
\end{equation}
(i.e. $k$ is the smallest integer such that $2k\ge n_E$). We then
use the corresponding projection $\projhat{E}{}$, see \eqref{eq:defProj0} and \eqref{eq:defProj},
to define 
\begin{equation}
\label{eq:defahE}
\ahE{\bs\sigma_h}{\bs\tau_h} = \scal[E]{\projhat{k-1}{E}\bs\sigma_h}{\projhat{k-1}{E}\bs\tau_h} \quad \forall \bs\sigma_h,\,\bs\tau_h \in \Vh[E].
\end{equation}
\begin{remark}\label{rm:dimension} 
We remark that, although a rigorous analysis is still missing for general polygons, the numerical tests (see Section \ref{subsec:convPoisson}) seem to suggest that the choice \eqref{k-choice} always leads to a stable scheme.
\end{remark}

\begin{remark}\label{rm:quadrature}
Given $\bftau \in \Vh[E]$, to compute $\projhat{k-1}{E}\bs\tau_h$ one has solve, from \eqref{eq:defProj} and integrating by parts:

\begin{equation}~\label{eq:defProj2}
\scal[E]{\projhat{k-1}{E}\bs\tau_h}{\nabla p} = -\scal[E]{\div \bs\tau_h}{p} + \int_{\partial E} (\bftau_h\cdot\bbn) p~\de, \qquad \quad \forall \; p\in \PolyH{k}{E},
\end{equation}
which is clearly computable, as $\div\bftau_h$ is computable and constant. Moreover, integrating by parts also the left-hand side and taking into account that the involved polynomials are harmonic, one realizes that the integral over $E$ can be computed as an integral over $\partial E$; therefore, only 1D quadrature rules are required to compute the left-hand side of \eqref{eq:defProj2}. Furthermore, since $\div\bftau_h$ is constant, the first term in the right-hand side requires only to evaluate the integral of a harmonic polynomial of degree at most $k$. Hence, the computation of $\projhat{k-1}{E}\bs\tau_h$ is not as cumbersome as it may appear at a first sight.       
\end{remark}

\paragraph{The local right-hand side term} We split the right-hand side term on each quadrilateral and we have  
\begin{equation*}
(f,v_h)_E= \int_{E}f v_h~\dEl.
\end{equation*}
Since $v_h\in\Uh(E)=\Poly{0}{E}$, we have that 
\begin{equation*}
(f,v_h)=\sum_{E\in\Mh}v_h \int_{E}f~\dEl,
\end{equation*}
which is computable via quadrature rules for polygonal domains, see for instance~\cite{Sommariva2007441}.
\subsection{The discrete scheme}\label{subsec:E2VEMdiscreteScheme}
Starting from the local spaces and local terms introduced in the previous sections, we can set the global self-stabilized problem.
More specifically, we introduce these two global approximation spaces, by gluing the local approximation
spaces, see~\eqref{eq:vemDiscreteSpace} and~\eqref{eq:polynomialDiscreteSpace}:
\begin{equation}\label{eq:globalVemSpace}
\Vh = \left\{ \bs{\tau}_h\in\sobhDiv{\Omega} : \bs{\tau}_{h|_{E}}\in \Vh[E],\quad \forall E\in\Mh \right\}
\end{equation}
and
\begin{equation}\label{eq:globalPolySpace}
\Uh=\left\{u_h\in U: u_{h|_{E}}\in \Uh(E),\quad  \forall E\in\Mh  \right\}.
\end{equation}
Now, given a local approximation of $a^E(\cdot,\cdot)$, see~\eqref{eq:defahE}, $\forall\, \bs\sigma_h,\,\bs\tau_h \in \Vh$ we set
\begin{equation}
\label{eq:defah}
\a[h]{\bs\sigma_h}{\bs\tau_h} := \sum_{E\in\Mh}\ahE{\bs\sigma_h}{\bs\tau_h}.
\end{equation}
We can state the discrete problem as: find $\left(\bs\sigma_h, u_h\right)\in\Vh\times\Uh$ such that
\begin{equation}
\label{eq:discrVarForm}
\begin{cases}
\ah{\bs\sigma_h}{\bs\tau_h} + \scal[\Omega]{\div\bs\tau_h}{u_h} = 0 \quad &\forall\bs{\tau}_{h}\in\Vh \\
\scal[\Omega]{\div\bs\sigma_h}{v_h} =\scal[\Omega]{f}{v_h}
\quad &\forall v_h \in \Uh
\end{cases}.
\end{equation}
In the next section we focus on the well-posedness of this discrete scheme, in the case of quadrilateral meshes,
which requires in particular the coercivity-on-the-kernel condition for the bilinear form $\ah{\cdot}{\cdot}$ (also called \emph{ellipticity-on-the-kernel condition}).

\section{Well-posedness in the quadrilateral case}
\label{sec:well-posedness}
\begin{figure}[]
		\centering
		\includegraphics[width=0.6\textwidth]{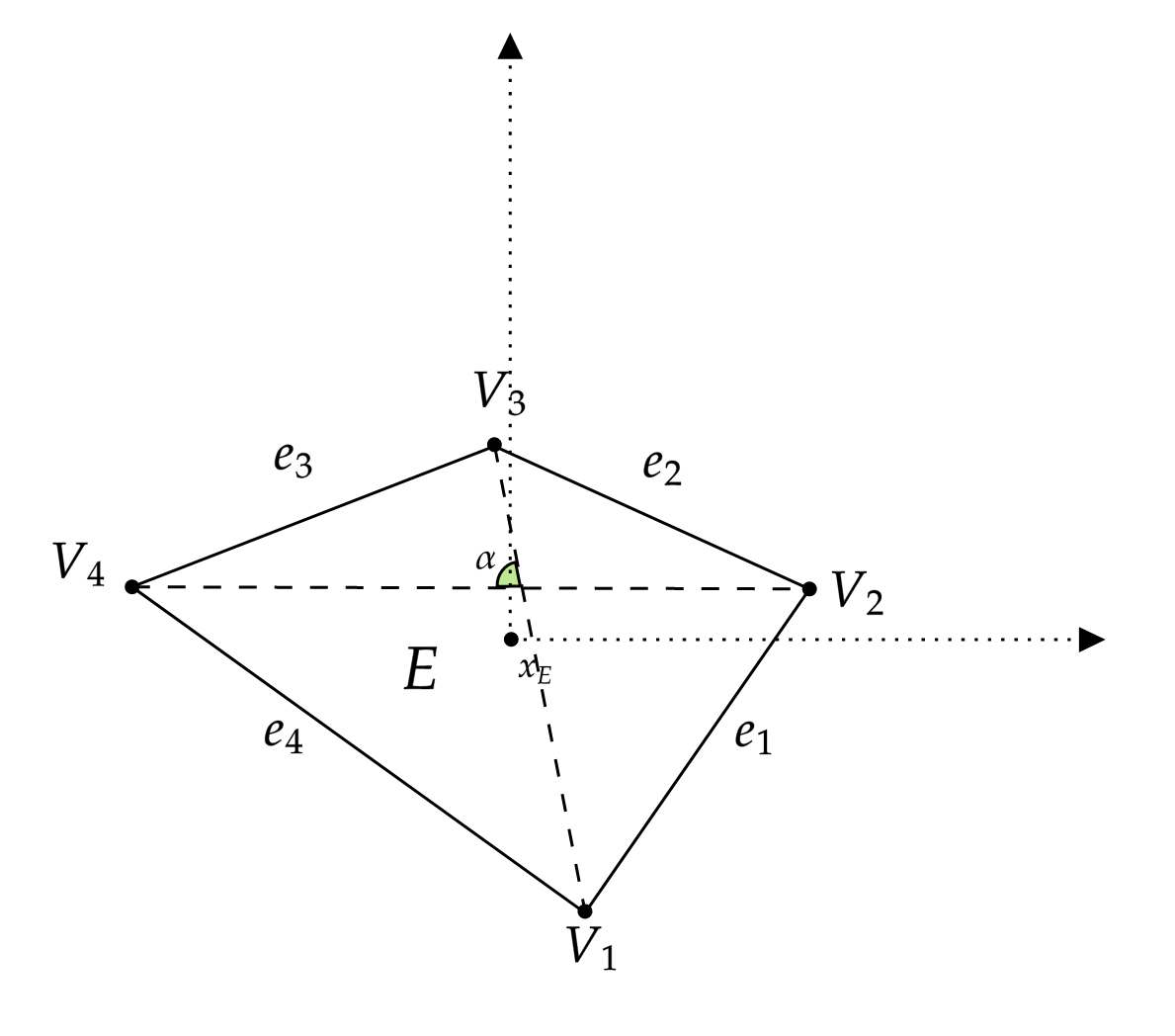}
		\caption{A general quadrilateral $E\in\Mh$}
		\label{Quadrilateral}
	\end{figure}

From now on we focus on the case where the mesh $\Mh$ is made up by quadrilaterals. This implies that we choose $k=2$, so that we use the local projection $\projhat{1}{E}$, see \eqref{k-choice}. Hence, we project onto the gradients of quadratic harmonic polynomials, a space of dimension $4$.
For each quadrilateral $E\in\Mh$, $V_i$ (for $i=1,\dots,4$) denote its vertices counterclockwise ordered and $e_i$ the edge connecting $V_i$ to $V_{i+1}$, where $V_5=V_1$ (see Figure.~\ref{Quadrilateral}). Let $\bbn_i$, be the unit normal vector of the edges $e_i$ for $i=1,\dots,4$.
This section is devoted to prove the well-posedness of the discrete problem stated by~\eqref{eq:discrVarForm}.

We introduce the following two useful spaces $\RT{E}$ and $\HG{E}$, and we prove some properties of their functions.
\begin{definition}[Raviart-Thomas space $\RT{E}$]
It is the space of the polynomial functions defined as follows:
\begin{equation}\label{eq:space_RT0}
\RT{E}:=
\left\{\bbr\in\lebldouble{E}\,:\, 
\bbr =
\begin{pmatrix}
c_1\\
c_2
\end{pmatrix}
+c_3
\begin{pmatrix}
x\\
y
\end{pmatrix},
\quad \mbox{s.t.} \quad c_1,c_2,c_3 \in \R
  \right\},
\end{equation}
whose dimension is equal to 3. 
\end{definition}
\begin{definition}[Hourglass space $\HG{E}$]
 Let $\hourglass\in\Vh[E]$ be the function such that
\begin{equation}\label{eq:fun_HG}
  \hourglass\cdot\bbn_{j} = \frac{(-1)^j}{|e_j|}
  \quad \forall j =1,\ldots,4,
\end{equation}
then we introduce the following one dimensional virtual space
\begin{equation}\label{eq:space_HG}
\HG{E}:= \mbox{span}\left( \hourglass \right).
\end{equation}
Using the divergence theorem, it is straightforward to see that a function $\bs\ttauh\in\HG{E}$ satisfies $\div\bs\ttauh=0$. 
\end{definition}
\begin{remark}
	We notice that the two spaces above are two subspaces of $\Vh(E)$.
\end{remark}

\begin{proposition}\label{pr:decomposition}
Let $\RT{E}$ be the space defined in \eqref{eq:space_RT0}
and let $\HG{E}$ be the space defined in \eqref{eq:space_HG},
then
\begin{equation}~\label{eq:decomposition_Sigma_h}
\Vh[E] = \RT{E} \oplus \HG{E}\,.
\end{equation}
Moreover, let us define the local divergence-free subspace:
\begin{equation}\label{eq:divfree_sp}
\Vhdf[E] = \left\{ {\bs\tau_h}\in \Vh[E]\ : \ {\div\bs\tau_h} = 0 \right\}.
\end{equation}
Then it holds
\begin{equation}~\label{eq:decomposition_Sigma_h_divfree}
\Vhdf[E] = (\Poly{0}{E})^2 \oplus \HG{E}\,
\end{equation}
and the decomposition is $L^2$-orthogonal.
\end{proposition}
\begin{proof}
Notice that, according to the dimension of $\RT{E}$ and $\HG{E}$, to get \eqref{eq:decomposition_Sigma_h}
we only have to prove that $\RT{E}\cap \HG{E}=\{\bs{0}\}$, that is $\hourglass \notin \RT{E}$.
By contradiction, we suppose that $\hourglass \in \RT{E}$. 
Notice that by definition of $\hourglass$, $\div \hourglass=0$,
hence $\hourglass\in (\Poly{0}{E})^2$. Take now $\bba = \nabla (\bba\cdot\bbx)$, where $\bba\in (\Poly{0}{E})^2$. We have, using integration by parts and \eqref{eq:fun_HG}:
\begin{equation}\label{eq:xinorm1}
(\hourglass, \bba)_E = (\hourglass,\nabla (\bba\cdot\bbx))_E = \int_{\partial E} (\hourglass\cdot \bbn_E)(\bba\cdot\bbx)= \sum_{i=1}^4 \int_{e_j} \frac{(-1)^j}{|e_j|} (\bba\cdot\bbx)\, , 
\end{equation}
for every $\bba\in (\Poly{0}{E})^2$.
An application of the trapezoidal rule gives
\begin{equation}\label{eq:xinorm2}
(\hourglass, \bba)_E   =  \bba \cdot \left(\frac{1}{2}\sum_{j=1}^4 (-1)^j ( V_j  + V_{j+1}) \right) = 0
\qquad \forall\, \bba\in (\Poly{0}{E})^2.
\end{equation}
Recalling that $\hourglass$ is constant, from \eqref{eq:xinorm2}
we infer $\hourglass = \bs{0}$, a contradiction since $\hourglass \ne 0$.
Furthermore, decomposition \eqref{eq:decomposition_Sigma_h_divfree} follows from a dimensional count, while the $L^2$-orthogonality is simply \eqref{eq:xinorm2}. 
\end{proof}

\begin{lemma}\label{lemma:HG_estimate}
	Let $E\in\Mh$ and let $\hourglass$ be the hourglass function defined on $E$ by \eqref{eq:fun_HG}. Then $\exists C_{\hourglass}>0$ independent of $h_E$ such that
	\begin{equation}
	\label{eq:HGfun_estim}
	\norm[0]{\hourglass} \leq C_{\hourglass} \,.
	\end{equation}
\end{lemma}
\begin{proof}
	Since $\hourglass\in\Vh[E]$, by \eqref{eq:vemDiscreteSpace}
	$\exists v \in\sobh{1}{E}$ such that $\hourglass = \nabla v$. It is
	clear that $v$ is defined up to a constant, so we choose $v$ such that
	$\int_E v = 0$. This implies that $\exists C > 0$ independent of $h_E$ such
	that
	\begin{equation}
	\label{eq:HGfun_estim:poincare-on-v}
	\norm[0]{v} \leq C h_E \norm[0]{\nabla v} = C h_E
	\norm[0]{\hourglass},
	\end{equation}
	by Poincar\'e's inequality. Moreover, since $\div\hourglass = 0$, it
	holds $\Delta v = 0$. Then, by Green's theorem and a Cauchy-Schwarz inequality
	we have
	\begin{equation}
	\label{eq:HGfun_estim:first-step}
	\norm[0]{\hourglass }^2 = \scal[E]{\hourglass }{\nabla v }
	= \scal[\partial E]{\hourglass\cdot \bbn}{ v} \leq
	\norm[0,\partial E]{\hourglass\cdot \bbn} \norm[0,\partial E]{v} \,.
	\end{equation}
	We can apply a standard trace inequality to the last norm and obtain, by
	exploiting also \eqref{eq:HGfun_estim:poincare-on-v},
	\begin{equation}
	\label{eq:HGfun_estim:traceineq-v}
	\norm[0,\partial E]{v} \leq h_E^{\frac12}\left(h_E^{-2}\norm[0]{v}^2
	+ \norm[0]{\nabla v}^2 \right)^{\frac12}
	\leq C h_E^{\frac12} \norm[0]{\hourglass } \,.
	\end{equation}
	On the other hand, an explicit computation exploiting the definition of
	$\hourglass $ given by \eqref{eq:fun_HG} yields
	\begin{equation}
	\label{eq:HGfun_estim:taudotn-estim}
	\norm[0,\partial E]{\hourglass \cdot\bbn}^2  = \sum_{j=1}^4 \int_{e_j}
	\left[
	\frac{(-1)^j}{|e_j|}
	\right]^2 = \sum_{j=1}^4 |e_j|^{-1} \leq 4\gamma^{-1} h_E^{-1} \,,
	\end{equation}
	where the last inequality is obtained by exploiting the mesh assumption~\ref{meshA_2}.  Using
	\eqref{eq:HGfun_estim:traceineq-v} and \eqref{eq:HGfun_estim:taudotn-estim} into
	\eqref{eq:HGfun_estim:first-step}, we get
	\begin{equation*}
		\norm[0]{\hourglass }^2
		\leq \norm[0,\partial E]{\hourglass\cdot \bbn} \norm[0,\partial E]{v}
		\leq 2\gamma^{-\frac12} h_E^{-\frac12} \cdot C h_E^{\frac12} \norm[0]{\hourglass }
		\leq C \norm[0]{\hourglass } \,,
	\end{equation*}
	which yields the thesis.
\end{proof}

\begin{lemma}\label{lem:infsupHG}
	Under the mesh assumptions \ref{meshA_1} and \ref{meshA_2}, for every $E\in\Mh$, there exists a positive constant $C_*$, independent of $h_E$, such that 
	\begin{equation}~\label{eq:infsupHGfunction}
			\norm[0]{\projhat{1}{E}\tilde{\bftau}_{h}} \geq C_*\norm[0]{\tilde{\bftau}_{h}} \qquad \forall\tilde{\bftau}_{h}\in\HG{E}.
	\end{equation}
\end{lemma}

\begin{proof}

Since $\HG{E}={\rm span}(\hourglass)$, it is sufficient to prove \eqref{eq:infsupHGfunction} for 
$\tilde{\bftau}_{h}=\hourglass$. 
Using the definition of the norm of the operator $\projhat{E}{}$
and~\eqref{eq:defProj}, 
we have
\begin{equation}~\label{eq:infsup_1}
		\norm[0]{\projhat{1}{E}\hourglass}
  = \sup_{\bbq
			\in\nabla\PolyH{2}{E}}
		\frac{\left(\projhat{1}{E}\hourglass,\bbq\right)}{\norm[0]{\bbq}}
		=\sup_{\bbq \in\nabla\PolyH{2}{E}}
		\frac{\left(\hourglass,\bbq\right)}{\norm[0]{\bbq}}.
\end{equation}

By Varignon's theorem~\cite{Coxeter1967}, for each element $E\in\Mh$, the quadrilateral $K_E$ whose vertices are the edge midpoints $M_j\;(j=1,\ldots,4)$ of $E$, is a parallelogram. With the usual abuse of notation that $V_5= V_1$ we have
\begin{equation*}
M_j = \frac{V_j+V_{j+1}}{2}\,,
\end{equation*}
 and the area of $K_E$ satisfies $\abs{K_E}=\frac{\abs{E}}{2}$.
Under the mesh assumptions~\ref{meshA_1} and~\ref{meshA_2}, it is not hard to show that the parallelogram is not degenerate, i.e. assumptions~\ref{meshA_1} and~\ref{meshA_2} hold for $K_E$ as well.
We now construct $p^*\in\PolyH{2}{E}$ such that
\begin{equation}\label{eq:altern}
	p^*\left(M_j\right)=(-1)^j, \quad \mbox{for each } j=1,\ldots,4.
\end{equation}
To this aim, it is useful to resort to complex numbers $z=x+i y$. Hence, up to a translation, we can identify $M_1$ as $0\in \C$; accordingly, we also set $M_2=z_1$, $M_4=z_2$ and $M_3 = z_1+z_2$. A direct computation shows that the complex-valued polynomial 
$$
q(z) = - 1 + 2 \frac{z_1+z_2}{z_1 z_2}\, z -  \frac{2}{z_1 z_2}\, z^2
$$
satisfies conditions \eqref{eq:altern} (with the above-mentioned identifications of $M_j$). 
We now set $p^*(x,y) = \mathrm{Re}(q(z))$, where $z = x + i y$ and $\mathrm{Re}(\cdot)$ denotes the real part. The real-valued polynomial $p^*$ is harmonic and satisfies conditions \eqref{eq:altern} as well.
Let $\bbp^\ast:=\nabla p^\ast$;
from~\eqref{eq:infsup_1} we get
\begin{equation}
\norm[0]{\projhat{1}{E}\hourglass}
	\geq  \frac{\left(\hourglass,\bbp^*\right)}{\norm[0]{\bbp^*}}.
\end{equation}
By an explicit computation using Cavalieri-Simpson's quadrature rule and \eqref{eq:altern}, we have that
  \begin{equation}\label{eq:integral_pstar}
	\begin{aligned}
		\left(\hourglass,\bbp^*\right) &= \int_E \hourglass\cdot \nabla p^*~\dEl
		=\int_{\partial E} (\hourglass\cdot\bbn) p^*~\de
		= \sum_{j=1}^4
		\frac{(-1)^j}{|e_j|}\int_{e_j}p^*~\de\\
		&=\sum_{j=1}^4 \frac{(-1)^j}{6}  \left(p^*\left(V_j\right)+4p^*\left(M_j\right)+p^*\left(V_{j+1}\right) \right)\\
	&= \sum_{j=1}^4  \frac{2}{3}(-1)^j p^*\left(M_j\right)
	=   \frac{8}{3}.	
	\end{aligned}
\end{equation}

We now notice that, due to assumptions~\ref{meshA_1} and~\ref{meshA_2},
there exists $C_{\bbp^\ast}>0$, independent of $h_E$, such that $\norm[0]{\bbp^*} = \norm[0]{\nabla p^*} \leq C_{\bbp^\ast}$.
Therefore, using Lemma~\ref{lemma:HG_estimate} we have

\begin{equation}
\norm[0]{\projhat{1}{E}\hourglass}\geq 
 \frac{ 8}{3 \norm[0]{\bbp^*}} 
 = \frac{ 8 \norm[0]{\hourglass}}{3 \norm[0]{\hourglass}\norm[0]{\bbp^*}} 
 \geq \frac{8}{3 C_{\hourglass} C_{\bbp^\ast}}\, \norm[0]{\hourglass}  \,.
\end{equation}
Then, \eqref{eq:infsupHGfunction} holds with $C_\ast = \frac{8}{3 C_{\hourglass}C_{\bbp^\ast}}$.
\end{proof}

\subsection{Continuity and coercivity of the local bilinear form $a_h^{E}(\cdot,\cdot)$}
In this section, applying the above preliminary results, in particular Lemma \ref{lem:infsupHG}, we prove the continuity and coercivity (on the divergence operator kernel) of the local bilinear form $a_h^{E}(\cdot,\cdot)$ in the $\lebl{}$-norm.
\begin{theorem}~\label{theo:continuityAndCoercivity}
	Under the mesh assumptions \ref{meshA_1} and \ref{meshA_2}, for every $E\in\Mh$, the discrete bilinear form $a_h^E(\cdot,\cdot)$, defined in~\eqref{eq:defahE}, is $\mathrm{L}^2$ continuous and coercive-on-the kernel, namely there exist two positive constants $\alpha_*$ and $\alpha^*$, independent of $h_E$, such that
	\begin{equation}~\label{eq:continuity}
		a_h^E(\bftau_h,\bfsigma_h) \leq \alpha^* \norm[0]{\bftau_h}\norm[0]{\bfsigma_h},\qquad	\forall \bftau_h, \bfsigma_h \in \Vh(E)
	\end{equation}
	and
	\begin{equation}~\label{eq:coercivity}
		a_h^E(\bftau_h,\bftau_h) \geq \alpha_* \norm[0]{\bftau_h}^2,\qquad	\forall \bftau_h\in \Vhdf(E)\, ,
	\end{equation}
where $\Vhdf$ is the divergence-free subspace defined in \eqref{eq:divfree_sp}. 
\end{theorem}
\begin{proof}
Fixed an element $E\in\Mh$, we first check the continuity. For every
$\bftau_h,\bfsigma_h\in\Sigma_h(E)$, applying the definition of $\projhat{E}{}$, its continuity and the Cauchy-Schwarz inequality, we obviously obtain
\begin{equation}
	\begin{aligned}
		a_h^E(\bftau_h, \bfsigma_h)
		=(\projhat{1}{E}\bftau_h , \projhat{1}{E}\bfsigma_h)  
		\leq \norm[0]{\bftau_h}\norm[0]{\bfsigma_h}.
	\end{aligned}
\end{equation}
Then~\eqref{eq:continuity} holds with $\alpha^*=1$.

Now, we prove the $\Vhdf$-coercivity of the bilinear form $a_h^E(\cdot,\cdot)$. From Proposition \ref{pr:decomposition}, we get that every $\bftau_h\in\Vhdf$ can be written by means of the orthogonal decomposition
$$
\bftau_h = \bftau_0 + \tilde{\bftau}_{h}\, ,
$$
where $\bftau_0\in(\Poly{0}{E})^2$ and $\tilde{\bftau}_{h}\in \HG{E}$.
Moreover, one has 
$$
||\bftau_h||^2_0 = ||\bftau_0||^2_0 +||\tilde{\bftau}_{h}||^2_0\, .
$$

Using Lemma~\ref{lem:infsupHG}, the definition of the projection operator~\eqref{eq:defProj} and noticing that $\projhat{1}{E}\bftau_{0} = \bftau_0$ , we have 
\begin{equation}
\begin{aligned}
		a_h^E(\bftau_h,\bftau_h)&=
		\left(\projhat{1}{E}\bftau_h,\projhat{1}{E}\bftau_h\right)=\left(\projhat{1}{E}\bftau_{0} +\projhat{1}{E}\tilde{\bftau}_{h}  ,\projhat{1}{E}\bftau_{0} +\projhat{1}{E}\tilde{\bftau}_{h}\right)\\
		 &= \left(\bftau_{0},\bftau_{0}\right)
		 +2\left(\bftau_{0},\tilde{\bftau}_{h}\right)
		 +\left(\projhat{1}{E}\tilde{\bftau}_{h},\projhat{1}{E}\tilde{\bftau}_{h}\right)
		 \\
          &= \left(\bftau_{0},\bftau_{0}\right)
		 +\left(\projhat{1}{E}\tilde{\bftau}_{h},\projhat{1}{E}\tilde{\bftau}_{h}\right)
		 \\
		&\geq 
			\left(\bftau_{0},\bftau_{0}\right)
		+C_*\left(\tilde{\bftau}_{h},\tilde{\bftau}_{h}\right)
		 \\
		&\geq \min\left\{1,C_*\right\}\left[  ||\bftau_0||^2_0 +||\tilde{\bftau}_{h}||^2_0
		\right]\\
		&= C_*\norm[0]{\bftau_{h}}^2,
	\end{aligned}
\end{equation}
which yields the thesis, with $\alpha_*=C_*$.
			\end{proof}
\subsection{Ellipticity-on-the-kernel condition and inf-sup condition}
In this section, we consider the two conditions, i.e. the coercivity of the bilinear form $a_h^{E}(\cdot,\cdot)$ on the kernel of the mixed term and the LBB inf-sup condition, that imply the well-posedness of the discrete problem \eqref{eq:discrVarForm}.

Let us introduce the discrete kernel space given by
\begin{equation}
	K_h :=\left\{ \bs{\tau}_h\in \Vh\colon \scal[]{\div \bs{\tau}_h}{v_h}=0 \; \forall \, v_h\in U_h 
	\right\}\,.
\end{equation}
Notice that $\forall \bs{\tau}_h\in\, K_h$ we have that $\div\bs{\tau}_h=0$, so that $\bftau_{h|E}\in \Vhdf$ and $\norm[\Sigma]{\bs{\tau}_h}=\norm[0]{\bs{\tau}_h}$.
Hence, applying the local coercivity property \eqref{eq:coercivity} stated in Theorem \ref{theo:continuityAndCoercivity} and the definition of the bilinear form $\a[h]{\cdot}{\cdot}$ \eqref{eq:defah}, we obtain that $\exists C_\ast>0$, independent of $h$, such that
\begin{equation}\label{eq:ellipt-kernel}
	\a[h]{\bftau_h}{\bftau_h} \geq C_* \norm[\Sigma]{\bftau_h}^2,\qquad	\forall \bftau_h\in K_h.
\end{equation}

Furthermore, the inf-sup condition, i.e.  
$\exists \beta>0$, independent of $h$, such that
\begin{equation}\label{eq:infsup}
	\inf_{v\in U_h} \limits \sup_{\bs{\tau}_h\in\Vh} \limits \frac{\scal[\Omega]{\div\bs{\tau}_h}{v}}{\norm[0]{v}\norm[\Sigma]{\bs{\tau}_h}} \geq \beta\,.
\end{equation}
is a consequence of the so-called \emph{Fortin's trick}, cf. \cite{BrezziBoffiFortin}, when the interpolation operator $\Ih$ of Section \ref{subsec:interpolation} is considered (see in particular \eqref{eq:divCompat}, \eqref{eq:interpCont} and \eqref{eq:interpEstimateDivTau} with $s=0$).   

\section{Error estimates in the quadrilateral case}
\label{sec:apriori} 
We prove optimal a priori error estimates for the method presented in this work, when the mesh is made up by quadrilaterals. We remark that the proof follows the usual guidelines for the VEM mixed schemes; however, we provide all the details, for the sake of completeness.  

\begin{theorem}
Let $\left(\bs \sigma, u\right)\in \left[\sobh{1}{\Omega}\right]^2\times\sobh[0]{1}{\Omega}$ and $f\in\sobh{1}{\Omega}$ be respectively solution and forcing term of \eqref{eq:cont-var-form}. Then $\exists C>0$, independent of $h$, such that the unique solution $\left(\bs\sigma_h, u_h\right)\in\Vh\times\Uh$ of \eqref{eq:discrVarForm} satisfies the following error estimates:
\begin{align} \label{eq:estimateSigma}
\norm[0]{\bs\sigma-\bs\sigma_h}&\leq C h \seminorm[1]{\bs\sigma}\,,
\\ \label{eq:estimateDivSigma}
\norm[0]{\div(\bs\sigma-\bs\sigma_h)}&\leq C h \seminorm[1]{f} \,,
\\ \label{eq:estimateU}
\norm[0]{u-u_h}&\leq Ch\left(\seminorm[1]{u}+\seminorm[1]{\bs\sigma}\right)\,.
\end{align}
\end{theorem}

\begin{proof}
In order to prove \eqref{eq:estimateSigma}, 
let $\bssI:=\Ih\bs\sigma\in\Vh$ be the interpolant of $\bs{\sigma}$ defined in Section \ref{subsec:interpolation}.
Then, applying the triangle inequality we obtain
\begin{equation}\label{eq:firstStepSigma}
\norm[0]{\bs\sigma-\bs\sigma_h}\leq 
\norm[0]{\bs\sigma-\bssI}
+
\norm[0]{\bssI-\bssh} \,.
\end{equation}
Let us focus on the term $\norm[0]{\bssI-\bssh}$. 
Notice that, applying the second equation of discrete problem \eqref{eq:discrVarForm} and the property of the interpolant \eqref{eq:propertyProjdivInterp},
we have for each $E\in\Mh$
\begin{equation}\label{eq:divIntmSol}
\div\bssh = -\proj{0}{E} f = \proj{0}{E} \div\bs{\sigma}= \div \bssI
\imply \div\left(\bssI-\bssh\right) = 0\, ,
\end{equation}
hence $(\bssI-\bssh)_{|E}\in \Vhdf$ for each $E\in\Mh$ (therefore $(\bssI-\bssh)\in K_h$).
Notice that applying this relation to the first equation of the discrete problem \eqref{eq:discrVarForm} and to the first equation of the continuous problem \eqref{eq:cont-var-form} we obtain that $\ah{\bssh}{\bssI-\bssh} =0 $ and $a(\bfsigma,\bssI-\bssh)=0$.
Then, since $\bssI-\bssh\in K_h$ we can apply Theorem \eqref{theo:continuityAndCoercivity}, in particular \eqref{eq:coercivity}, and obtain the estimate
\begin{equation}
\begin{aligned}
\alpha_\ast & \norm[0]{\bssI-\bssh}^2 \leq \ah{\bssI-\bssh}{\bssI-\bssh}\\
&= \ah{\bssI}{\bssI-\bssh}
\\
&= \sum_{E\in\Mh} \left(\ahE{\bssI -\projhat{1}{E} \bs{\sigma} }{\bssI-\bssh} + \ahE{\projhat{1}{E} \bs{\sigma} }{\bssI-\bssh} \right) \\
&= \sum_{E\in\Mh} \left(\ahE{\bssI -\projhat{1}{E} \bs{\sigma} }{\bssI-\bssh} + \aE{\projhat{1}{E} \bs{\sigma} }{\bssI-\bssh}
\right) \\
&= \sum_{E\in\Mh} \left( \ahE{\bssI -\projhat{1}{E} \bs{\sigma} }{\bssI-\bssh} + \aE{\projhat{1}{E} \bs{\sigma} -\bs{\sigma}}{\bssI-\bssh} \right)
\,,
\end{aligned}
\end{equation}
where the projector $\projhat{1}{E}$ is defined by the orthogonality condition \eqref{eq:defProj} and satisfies, for each $E\in\Mh$, $\ahE{\projhat{1}{E}\bs{\sigma}}{\bs{\tau}_h}=\aE{\projhat{1}{E}\bs{\sigma}}{\bs{\tau}_h}$ $\forall \bs{\tau}_h\in\Vh$.
We now notice that, since $\bs\sigma = \nabla u$ and $\projhat{1}{E}$ projects onto the space $\nabla \PolyH{2}{E}$, it holds
\begin{equation}\label{eq:gradapprox}
\norm[0,E]{\bs{\sigma}-\projhat{1}{E} \bs{\sigma}} = \norm[0,E]{\nabla u -\projhat{1}{E} (\nabla u)}
=\inf_{p\in \PolyH{2}{E}} | u - p|_{1,E}\le C h_E |\bs\sigma|_{1,E},
\end{equation}
where the last estimate follows from the standard approximation theory, see \cite{BrambleHilbert,DupontScott,BrennerScott}. 
Then, by the continuity of $\ahE{\cdot}{\cdot}$ and $\aE{\cdot}{\cdot}$, applying estimates  \eqref{eq:interpEstimateTau} and \eqref{eq:gradapprox},  we obtain
\begin{equation}
\begin{aligned}
\norm[0]{\bssI-\bssh} &\leq C \sum_{E\in\Mh}\left( \norm[0,E]{\bssI -\projhat{1}{E} \bs{\sigma}} +\norm[0,E]{\bs{\sigma}-\projhat{1}{E} \bs{\sigma}}\right)
\\
&\leq C \left( \norm[0]{\bs{\sigma}-\bssI} +\sum_{E\in\Mh}\norm[0,E]{\bs{\sigma}-\projhat{1}{E} \bs{\sigma}}\right) 
\\
& \leq Ch\seminorm[1]{\bs{\sigma}}\,.
\end{aligned}
\end{equation}
Applying this relation and the interpolation estimate \eqref{eq:interpEstimateTau} to \eqref{eq:firstStepSigma},estimate \eqref{eq:estimateSigma} is proved.
Moreover, to prove \eqref{eq:estimateDivSigma} we apply \eqref{eq:divIntmSol}, the interpolation estimate \eqref{eq:interpEstimateDivTau} and the equation $\div\bs\sigma = -f$, to obtain
\begin{equation}
\norm[0]{\div\bs\sigma-\div\bs\sigma_h}
\leq \norm[0]{\div\bs\sigma-\div\bs\sigma_I} 
\leq C h \seminorm[1]{f}\,.
\end{equation}
Finally, we have to prove \eqref{eq:estimateU}. 
Let $u_I := \proj{0}{h} u\in U_h$.
Notice that by its definition $u_I$ satisfies $\scal[\Omega]{u-u_I}{\div\bs\tau_h}=0$ for each $\bs\tau_h\in\Vh$.
By triangle inequality we have
\begin{equation}\label{eq:step1-estimateU}
\norm[0]{u-u_h} \leq \norm[0]{u-u_I} +\norm[0]{u_I-u_h}\,.
\end{equation}
First, let us consider the term $\norm[0]{u_I-u_h}$. Since $u_I-u_h\in U_h$, according to the inf-sup condition \eqref{eq:infsup}, there exists $\bs\tau_h^\star\in\Vh$ be such that $\div\bs\tau_h^\star =u_I-u_h$ and 
\begin{equation}\label{eq:normRelDiv}
\norm[0]{\bs\tau_h^\star} \leq \frac{1}{\beta} \norm[0]{u_I-u_h}\,.
\end{equation}
Then, applying the continuous problem \eqref{eq:cont-var-form}, the discrete one \eqref{eq:discrVarForm} and adding and subtracting $\projhat{E}{} \bs{\sigma}$, we obtain
\begin{equation}\label{eq:estimateUIUh}
\begin{aligned}
\norm[0]{u_I-u_h}^2 
&= \scal[\Omega]{u_I-u_h}{\div\bs\tau_h^\star} \\
&= \scal[\Omega]{u-u_h}{\div\bs\tau_h^\star} \\
& = a(\bs{\sigma},\bs\tau_h^\star)- \ah{\bs\sigma_h}{\bs\tau_h^\star}
\\
& = \sum_{E\in\Mh} \aE{\bs\sigma-\projhat{1}{E} \bs{\sigma}}{\bs\tau_h^\star} +\ahE{\projhat{1}{E} \bs{\sigma}-\bssh}{\bs\tau_h^\star}
\\
&\leq C\sum_{E\in\Mh}\left(
\norm[0,E]{\bs\sigma-\projhat{1}{E} \bs{\sigma}} +\norm[0,E]{ \bs{\sigma}-\bssh}
\right)\norm[0]{u_I-u_h}
\end{aligned}
\end{equation}
where in the last step we exploit the continuity of the bilinear forms together with \eqref{eq:normRelDiv}.
Finally, applying \eqref{eq:estimateUIUh} to \eqref{eq:step1-estimateU}, the interpolation estimate \eqref{eq:interpEstimateU} and the error estimate of $\bs\sigma$ \eqref{eq:estimateSigma} already proved, we obtain
\begin{equation}
\begin{aligned}
\norm[0]{u-u_h} 
&\leq \norm[0]{u-u_I} + C \left(\sum_{E\in\Mh} \norm[0,E]{\bs\sigma-\projhat{1}{E} \bs{\sigma}}+
\norm[0]{\bs\sigma-\bssh}\right)
\\
&\leq C h \left(\seminorm[1]{u}+ \seminorm[1]{\bs\sigma}\right)\,.
\end{aligned}
\end{equation}
\end{proof}
\section{Numerical Tests}\label{sec:numtests}
\subsection{Convergence tests}\label{subsec:convPoisson}

\begin{figure} 
  \centering
  \begin{subfigure}[b]{.24\textwidth}
    \includegraphics[width=\linewidth]{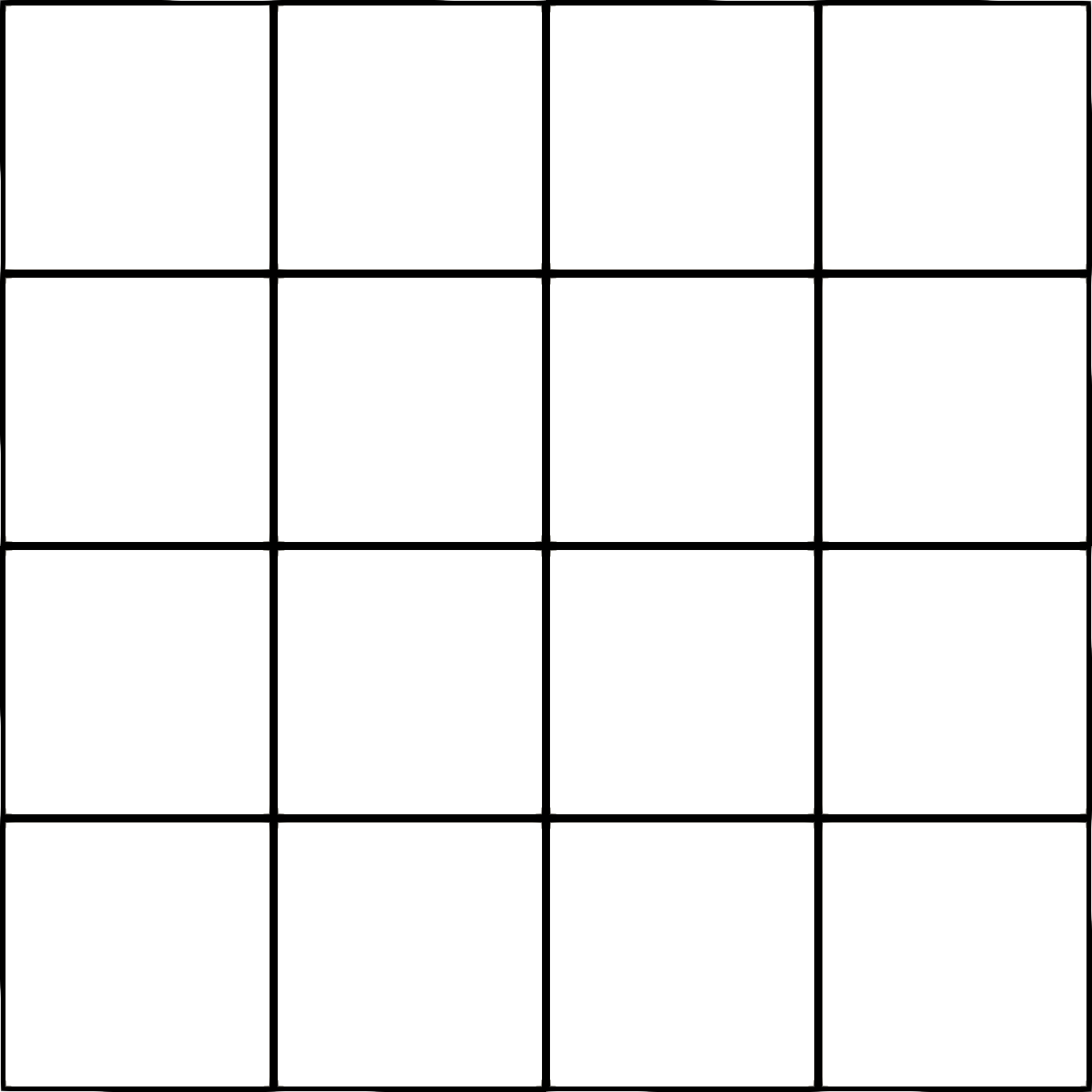}
    \caption{\meshtag{Cartesian}}\label{subfig:cartesian}
  \end{subfigure}
  \begin{subfigure}[b]{.24\textwidth}
    \includegraphics[width=\linewidth]{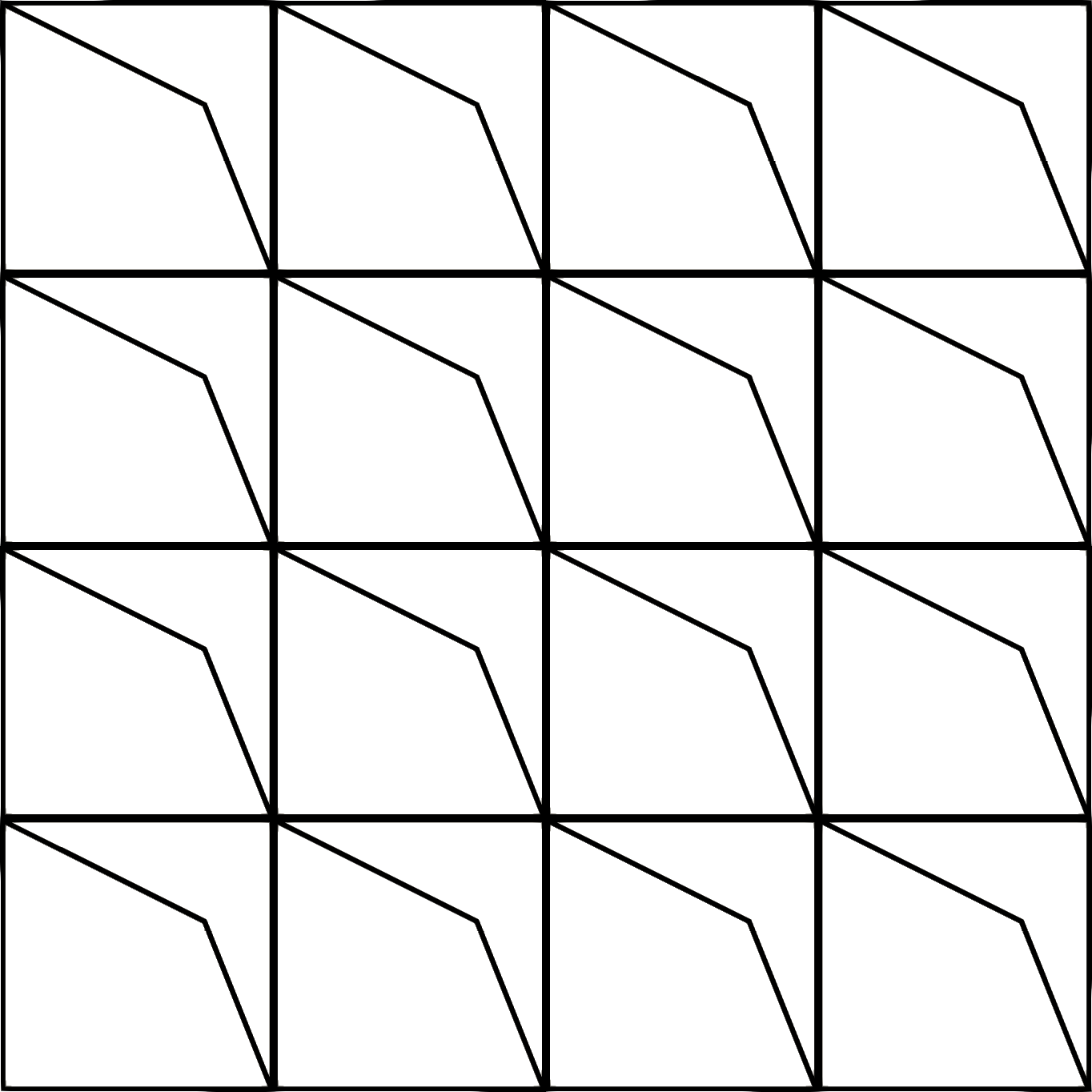}
    \caption{\meshtag{ConvexConcave}}\label{subfig:convexconcave}
  \end{subfigure}    
  \begin{subfigure}[b]{.24\textwidth}
    \includegraphics[width=\linewidth]{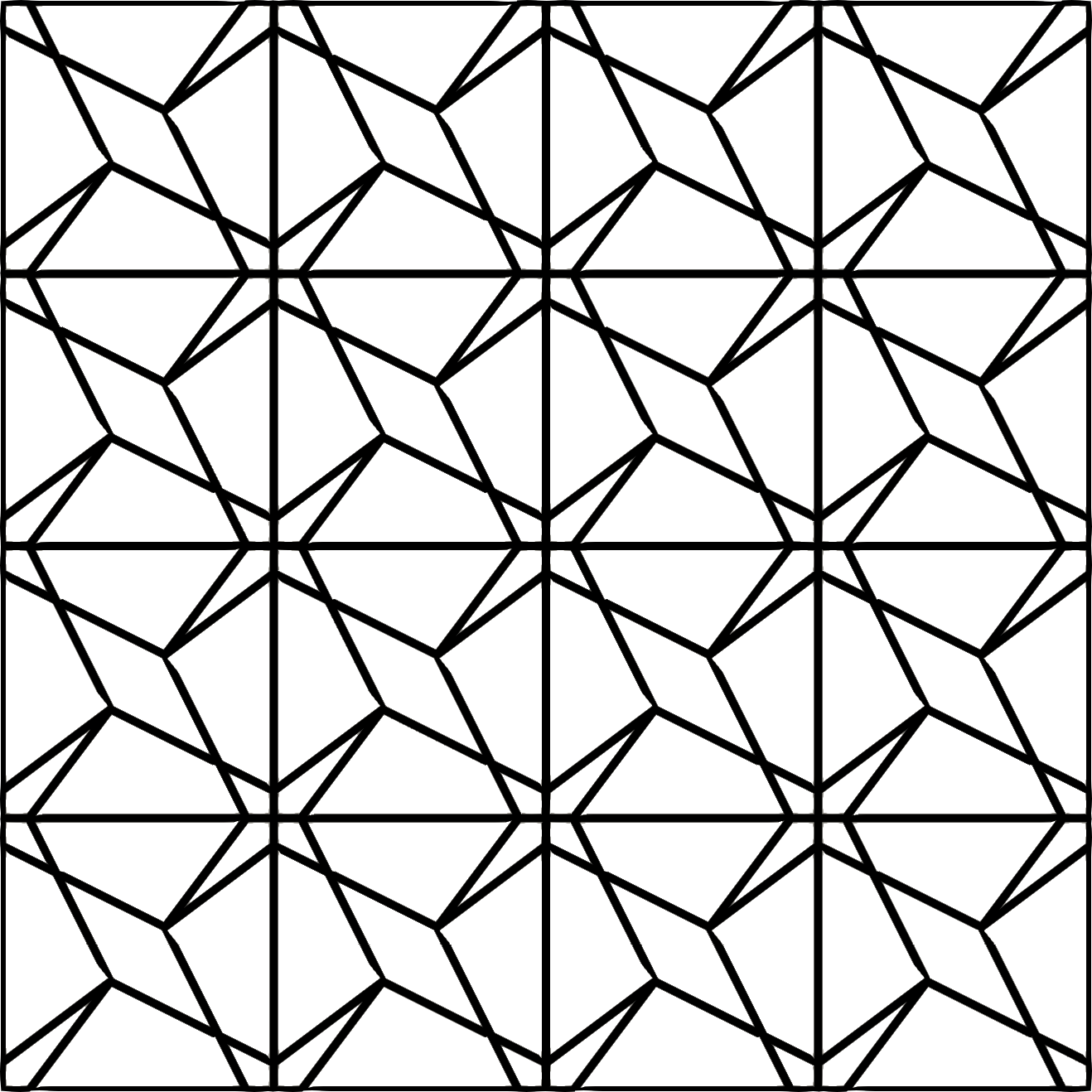}
    \caption{\meshtag{Distorted}}\label{subfig:distorted}
  \end{subfigure} 
    \begin{subfigure}[b]{.24\textwidth}
    \includegraphics[width=\linewidth]{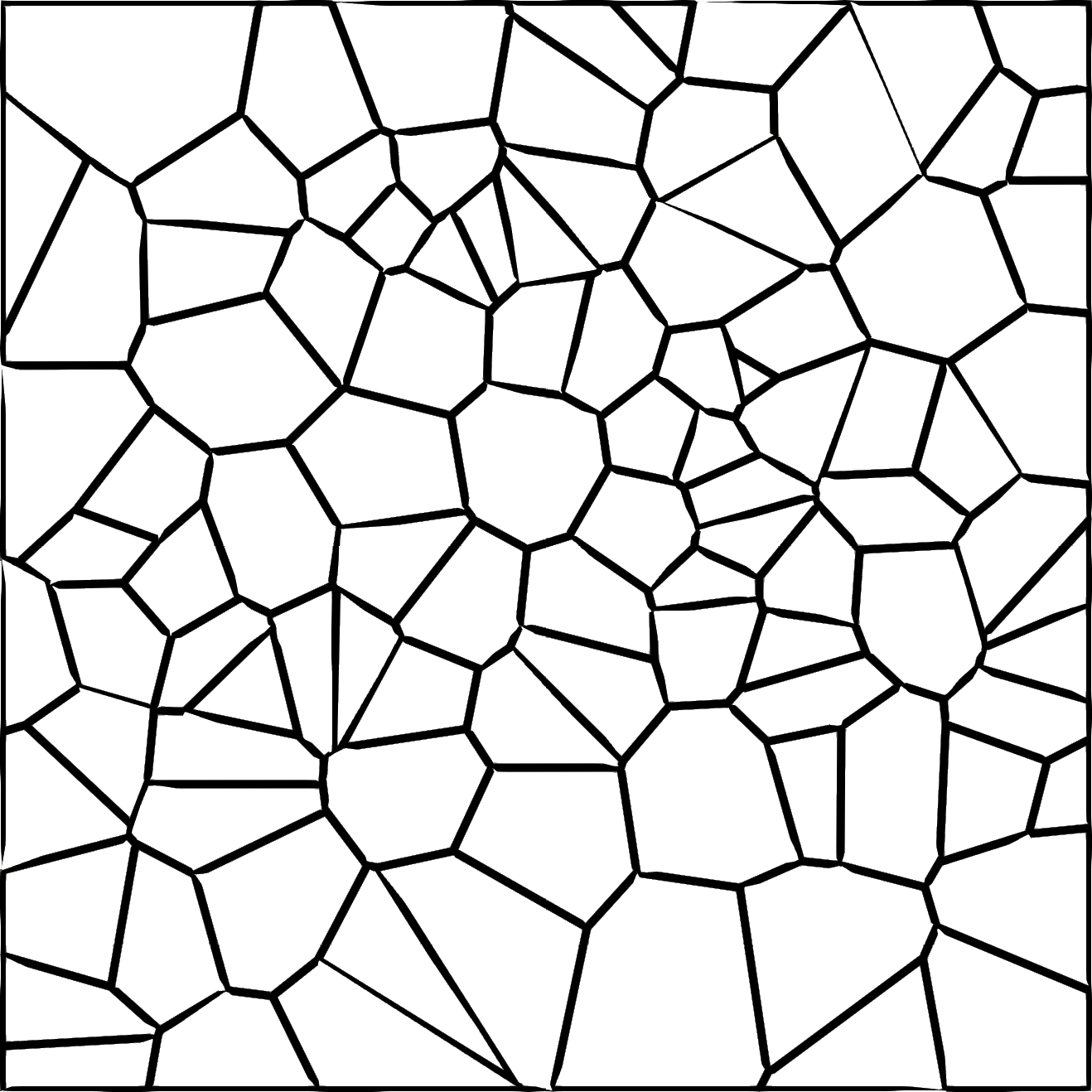} %
    \caption{\meshtag{Random}}\label{subfig:random}
  \end{subfigure} 
  \caption{Meshes.}
  \label{fig:meshes}
  \end{figure}

\begin{figure} 
  \centering
  \begin{subfigure}[b]{.32\linewidth}
    \includegraphics[width=\linewidth]{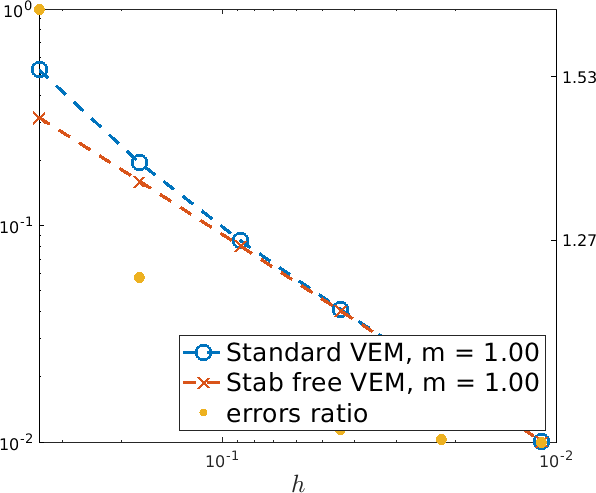}
    \caption{\meshtag{Cartesian}, $\mathrm{err}_u$}
  \end{subfigure}
  \begin{subfigure}[b]{.32\linewidth}
 		\includegraphics[width=\linewidth]{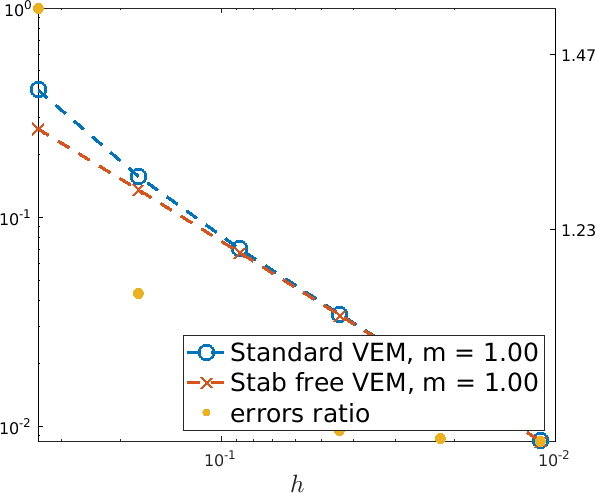}
 		\caption{\meshtag{ConvexConcave}, $\mathrm{err}_u$}
  \end{subfigure}
    \begin{subfigure}[b]{.32\linewidth}
 		\includegraphics[width=\linewidth]{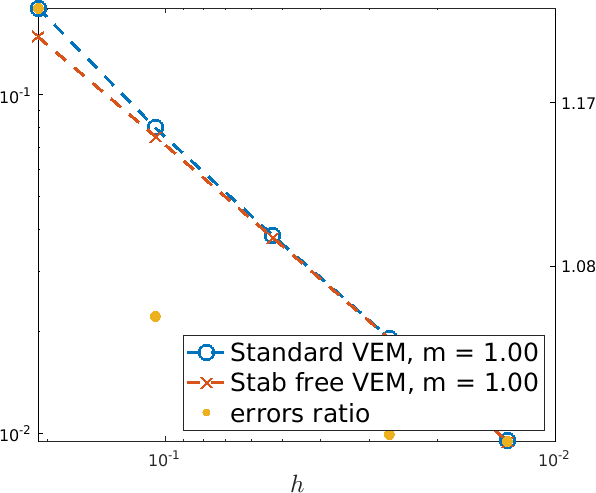}
 		\caption{\meshtag{Distorted}, $\mathrm{err}_u$}
  \end{subfigure}
  
    \begin{subfigure}[b]{.32\linewidth}
    \includegraphics[width=\linewidth]{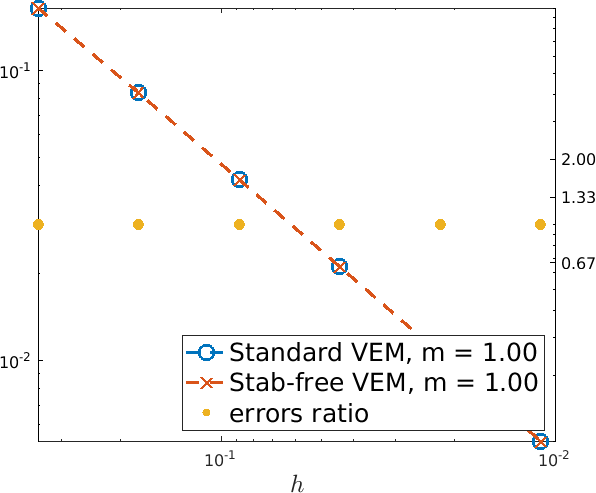}
    \caption{\meshtag{Cartesian}, $\mathrm{err}_{\div}$}
  \end{subfigure}
  \begin{subfigure}[b]{.32\linewidth}
    \includegraphics[width=\linewidth]{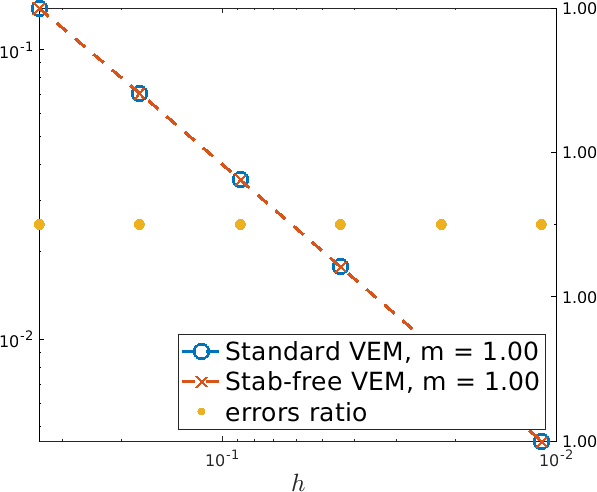}
    \caption{\meshtag{ConvexConcave}, $\mathrm{err}_{\div}$}
  \end{subfigure}
    \begin{subfigure}[b]{.32\linewidth}
    \includegraphics[width=\linewidth]{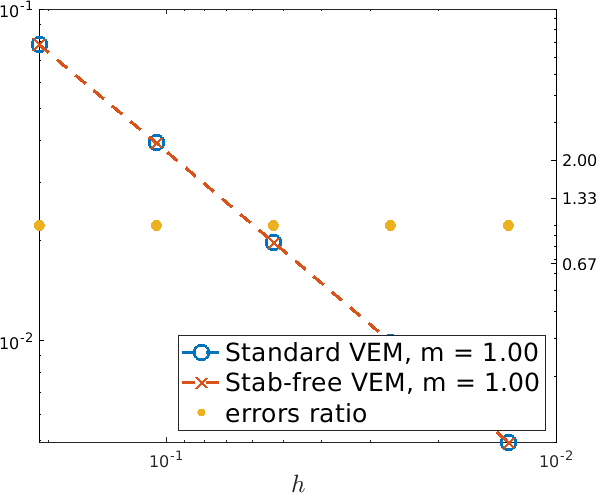}
    \caption{\meshtag{Distorted}, $\mathrm{err}_{\div}$}
  \end{subfigure}
  
\begin{subfigure}[b]{.32\linewidth}
	\includegraphics[width=\linewidth]{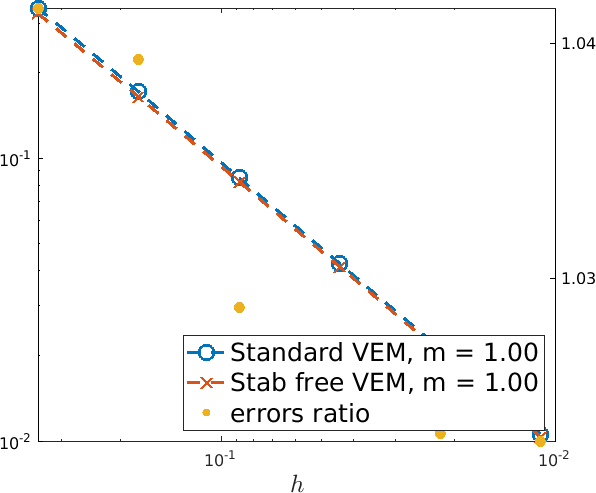}
	\caption{\meshtag{Cartesian}, $\mathrm{err}_{\bs\sigma}$}
\end{subfigure}
  \begin{subfigure}[b]{.32\linewidth}
    \includegraphics[width=\linewidth]{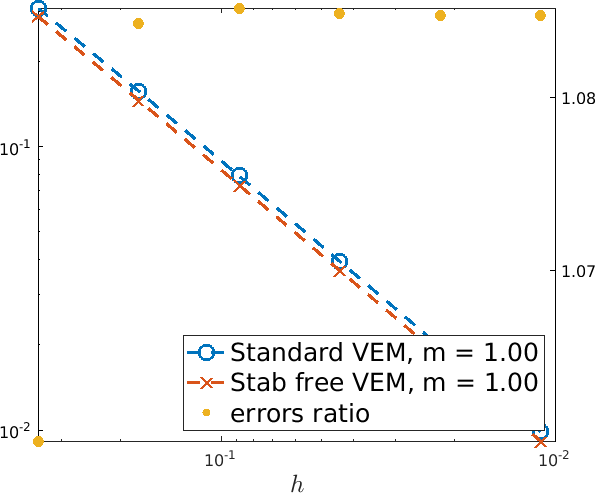}
    \caption{\meshtag{ConvexConcave}, $\mathrm{err}_{\bs\sigma}$}
  \end{subfigure}  
    \begin{subfigure}[b]{.32\linewidth}
    \includegraphics[width=\linewidth]{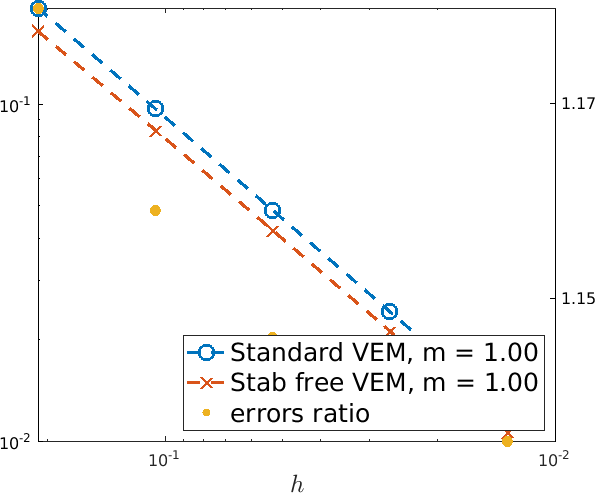}
    \caption{\meshtag{Distorted}, $\mathrm{err}_{\bs\sigma}$}
  \end{subfigure}  
  
\begin{subfigure}[b]{.32\linewidth}
	\includegraphics[width=\linewidth]{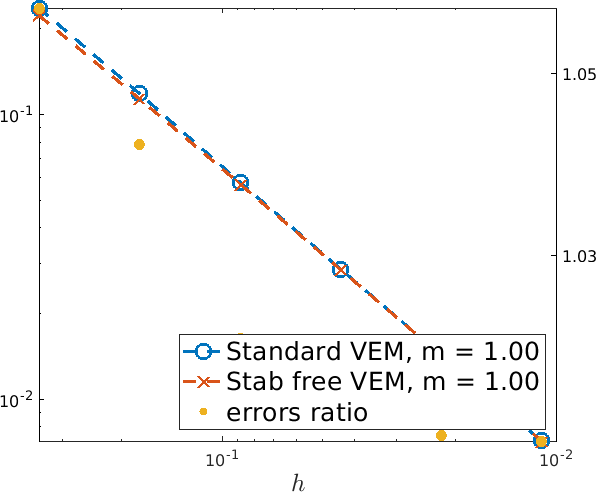}
	\caption{\meshtag{Cartesian}, $\mathrm{err}_{\bs\sigma\cdot\bs{n}}$}
\end{subfigure}
  \begin{subfigure}[b]{.32\linewidth}
    \includegraphics[width=\linewidth]{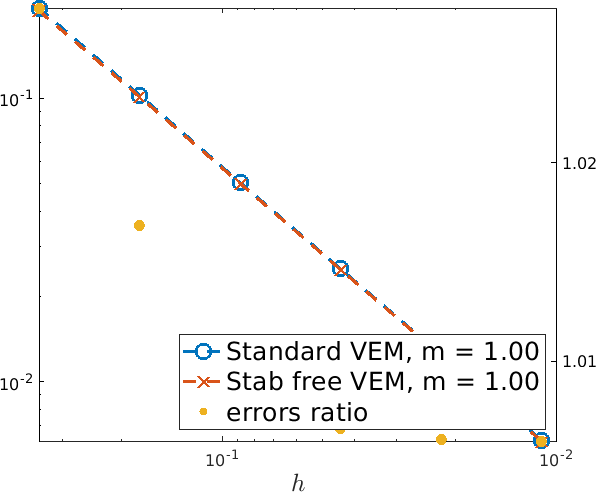}
    \caption{\meshtag{ConvexConcave}, $\mathrm{err}_{\bs\sigma\cdot\bs{n}}$}
  \end{subfigure}  
    \begin{subfigure}[b]{.32\linewidth}
    \includegraphics[width=\linewidth]{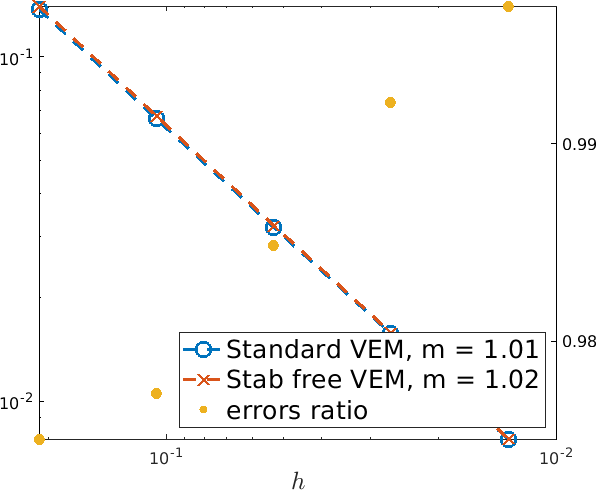}
    \caption{\meshtag{Distorted}, $\mathrm{err}_{\bs\sigma\cdot\bs{n}}$}
  \end{subfigure} 
  
  \caption{Convergence curves on quadrilateral meshes. The left vertical axis refers to the values of the errors (\textit{dotted lines}). The right vertical axis refers to the ratio between the error made by the standard VEM method and the error of the proposed method (\textit{orange dots}).}
  \label{fig:poissonConvergenceQuadrilateral}
\end{figure}

\begin{figure} 
  \centering
  \begin{subfigure}[b]{.32\linewidth}
    \includegraphics[width=\linewidth]{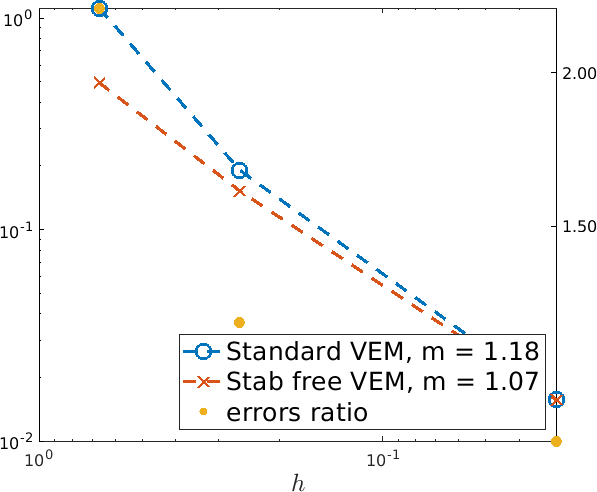}
    \caption{$\mathrm{err}_u$}
  \end{subfigure}

  \begin{subfigure}[b]{.32\linewidth}
    \includegraphics[width=\linewidth]{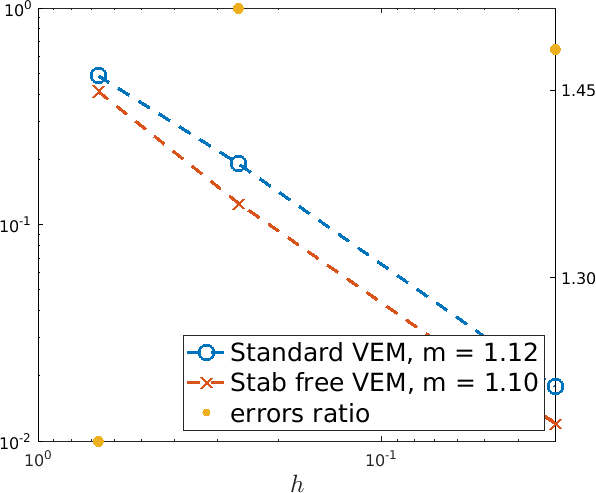}
    \caption{$\mathrm{err}_{\bs\sigma}$}
  \end{subfigure}  
  \begin{subfigure}[b]{.32\linewidth}
    \includegraphics[width=\linewidth]{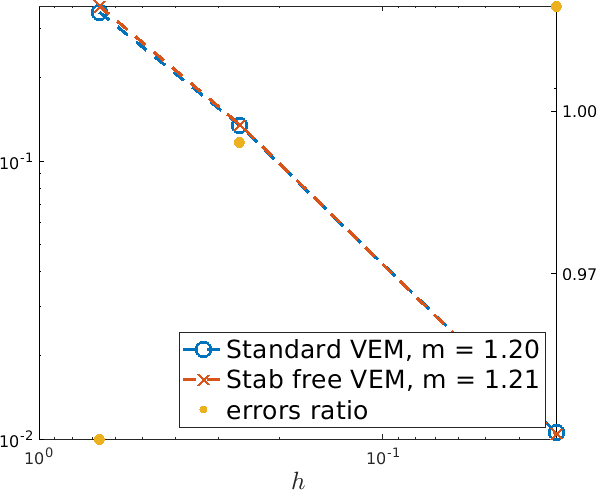}
    \caption{$\mathrm{err}_{\bs\sigma\cdot\bs{n}}$}
  \end{subfigure}
  \caption{Convergence curves on \meshtag{Random} mesh. The left vertical axis refers to the values of the errors (\textit{dotted lines}). The right vertical axis refers to the ratio between the error made by the standard VEM method and the error of the proposed method (\textit{orange dots}).}
  \label{fig:convPoissonRandom}
\end{figure}

In this section, we numerically assess the behaviour of our scheme with respect
to mesh refinement. We consider $\Omega = (0,1)^2$ and solve Problem
\eqref{eq:discrVarForm} choosing $f$ such that
\begin{align*}
  u(x,y) &=  x(1-x)y(1-y) \,,
  \\
  \bs\sigma(x,y) &= \nabla u (x,y) =
                   \begin{pmatrix}
                     (1-2x)y(1-y) \\ x(1-x)(1-2y)
                   \end{pmatrix}
  \,.
\end{align*}
First, we consider the four
families of meshes depicted in Figure
\ref{fig:meshes}. We assess the method behaviour by computing the following
relative errors:
\begin{align*}
  \mathrm{err}_u &= \frac{1}{\norm[0]{u}}
                   \left(
                   \sum_{E\in\Mh}\norm[0,E]{u - u_h}^2
                   \right)^{\frac12} \,,
  \\
  \mathrm{err}_{\div} &= \frac{1}{\norm[0]{\div\bs\sigma}} \sum_{E\in\Mh}
                        \left(
                        \norm[0,E]{\div\bs\sigma - \div\bs\sigma_h}^2
                        \right)^{\frac12} \,,
  \\
  \mathrm{err}_{\bs\sigma}
                 &= \frac{1}{\norm[0]{\bs\sigma}}\left(
                   \sum_{E\in\Mh}\norm[0,E]{\bs\sigma - \projhat{k-1}{E}\bs\sigma_h}^2
                   \right)^{\frac12}  \,,
        \\
            \mathrm{err}_{\bs\sigma\cdot\bs n}
                 &=\frac{\left(\sum_{e\in\mathcal{E}_h} h_e \norm[0,e]{(\bs\sigma-\bs\sigma_h)\cdot\bs n^e } \right)^{\frac12} }{\left(\sum_{e\in\mathcal{E}_h} h_e \norm[0,e]{\bs\sigma\cdot\bs n^e }\right)^{\frac12}}
                   \,,
\end{align*}
where $\mathcal{E}_h$ denotes the set all edges of $\Mh$.
We also
solve the test problem with the standard VEM method \cite{Beirao2016}. 
We recall that for this latter method, the local discrete bilinear form $a_h(\cdot,\cdot)$ is given by
\begin{equation}
\label{eq:standardVEM}
\ahE{\bs\sigma_h}{\bs\tau_h} = \scal[E]{\proj{0}{E}\bs\sigma_h}{\proj{0}{E}\bs\tau_h} + s^E\left( (I-\proj{0}{E} )\bs\sigma_h , (I-\proj{0}{E} )\bs\tau_h \right),
\end{equation}
where $s^E(\cdot,\cdot)$ is the local stabilization term. In matrix form, the stabilization term we choose is given by

\begin{equation}
\label{eq:stabmatrix}
\mathbf{S} = (\mathbf{I}- \mathbf{\Pi}^0)^T \mathbf{D}\, (\mathbf{I} - \mathbf{\Pi}^0 ),
\end{equation}  
where the matrix $\mathbf{\Pi}^0$ represents the projection onto the constant vector functions. Moreover, $\mathbf{D}$ is a diagonal matrix defined as 
\begin{equation}
	\mathbf{D}_{ii} = \max\left( h_E |e_i|, 
	\scal[E]{\proj{0}{E}\bs\varphi_i}{\proj{0}{E}\bs\varphi_i}
	\right)\,.
\end{equation}
Above, the functions $\bs\varphi_i$ denote the elements of the Lagrangian basis corresponding to the local degrees of freedom \eqref{eq:defDofSigma}. 
This choice is known as D-recipe stabilization, and it is inspired by the numerical assessment in \cite{d-recipe1, d-recipe2}.
We also notice that in computing the error $\mathrm{err}_{\bs\sigma}$, for the standard VEM we use the projection onto constants. 

In Figure \ref{fig:poissonConvergenceQuadrilateral}, we consider the quadrilateral meshes of Figures \ref{subfig:cartesian}, \ref{subfig:convexconcave} and \ref{subfig:distorted}, respectively named \meshtag{Cartesian}, \meshtag{ConvexConcave} and \meshtag{Distorted}. 
The
computed errors obtained by the two methods are compared with respect to the
maximum diameter of the mesh, denoted by $h$, and the asymptotical convergence rates are
reported in the legend. 
The results show that the two methods behave
equivalently on all meshes with respect to all the computed errors. We also remark that on all the meshes, both methods return exactly the same results for $\mathrm{err}_{\div}$. This is not surprising, since for all the meshes and both methods, from the second equation of \eqref{eq:discrVarForm} we get
$\div\bssh = -\proj{0}{E} f$, while $\div\bfsigma = - f$. Hence $\mathrm{err}_{\div}$ is always the $L^2$ error when the load term $f$ is approximated by piecewise constant functions. Accordingly, from now on we will not display that error quantity. 

In Figure \ref{subfig:random} we consider the family of meshes named \meshtag{Random}, i.e. polygonal meshes obtained using \textit{Polymesher} \cite{Polymesher}, whose elements are not only quadrilaterals. 
On each polygon, we construct the local bilinear form $\ahE{\cdot}{\cdot}$ \eqref{eq:defahE} choosing $k$ as in \eqref{k-choice} (see Remark \ref{rm:dimension}). 
As we can see in Figure \ref{fig:convPoissonRandom}, the proposed method is stable and exhibits the expected convergence rates.

\subsection{Comparison with standard VEM on an anisotropic refinement test}
\begin{figure} 
  \centering
    \includegraphics[width=.33\linewidth]{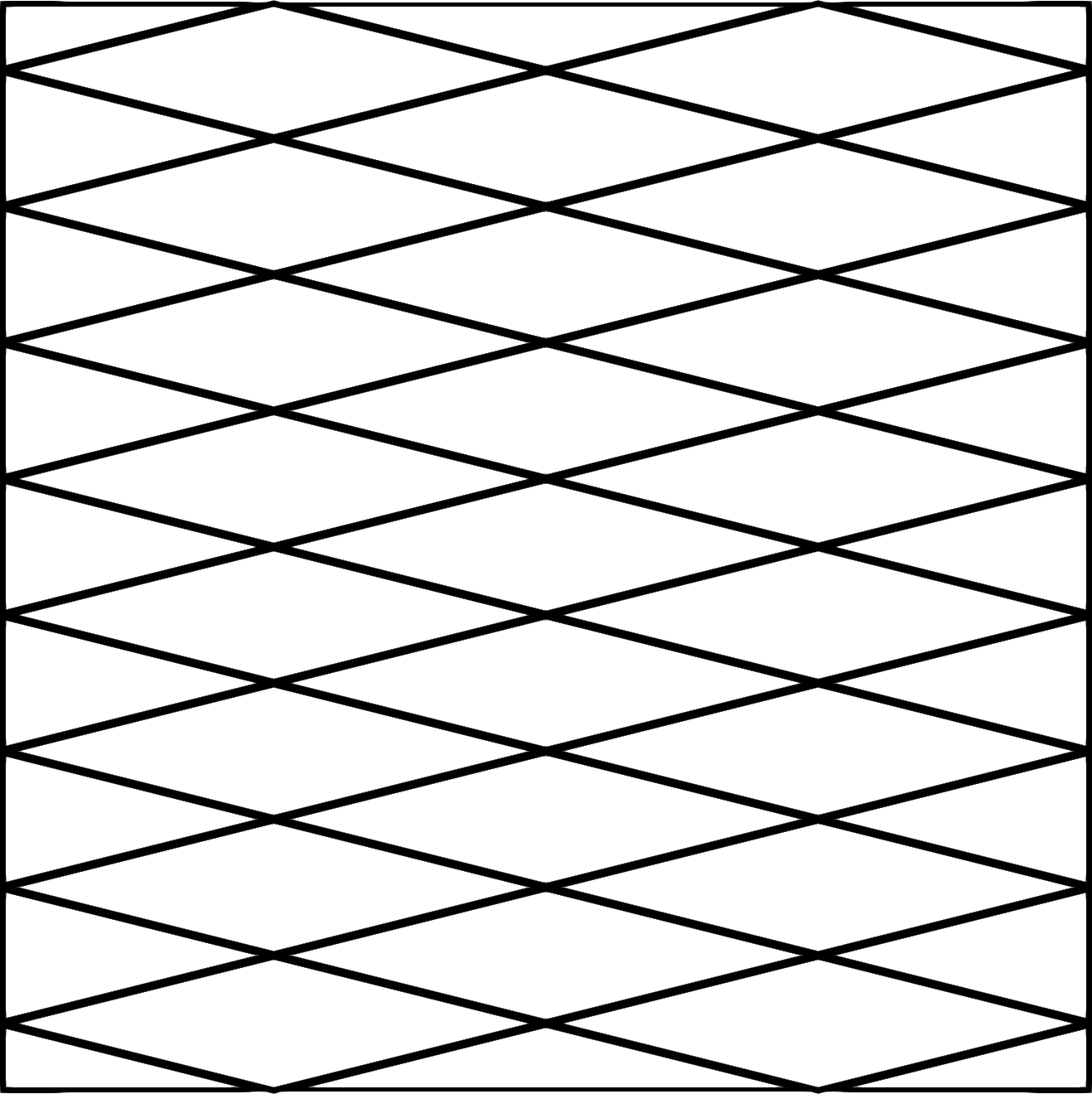} %
    \caption{\meshtag{Rhomboidal} mesh} 
  \label{fig:RhombiMesh}
  \end{figure}
  
\begin{figure} 
  \centering
  \begin{subfigure}[b]{.32\linewidth}
    \includegraphics[width=\linewidth]{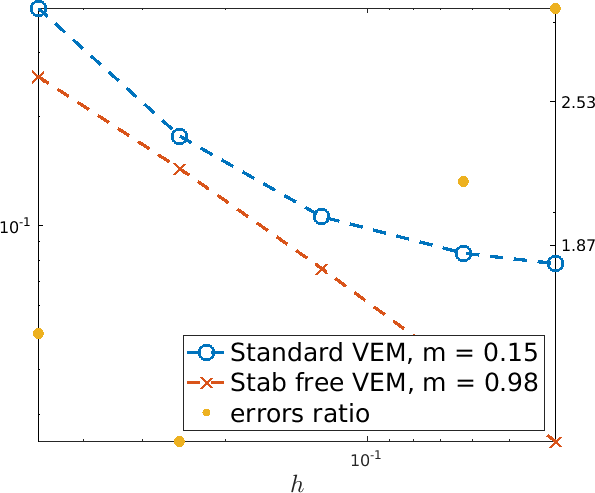}
    \caption{$\mathrm{err}_u$}
  \end{subfigure}
  \begin{subfigure}[b]{.32\linewidth}
    \includegraphics[width=\linewidth]{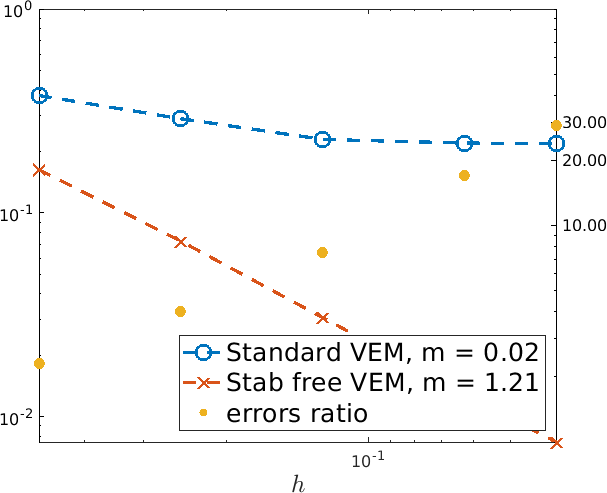}
    \caption{$\mathrm{err}_{\bs\sigma}$}
  \end{subfigure}  
  \begin{subfigure}[b]{.32\linewidth}
    \includegraphics[width=\linewidth]{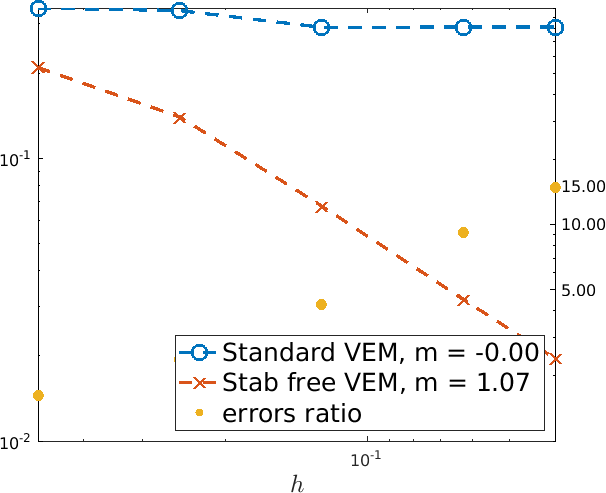}
    \caption{$\mathrm{err}_{\bs\sigma\cdot\bs{n}}$}
  \end{subfigure}
  \caption{Convergence curves on \meshtag{Rhomboidal} mesh. The left vertical axis refers to the values of the errors (\textit{dotted lines}). The right vertical axis refers to the ratio between the error made by the standard VEM method and the error of the proposed method (\textit{orange dots}).}
  \label{fig:convPoissonRhombi}
\end{figure}

In this section, we consider the problem presented in the previous section with a \meshtag{Rhomboidal} mesh, as depicted in Figure \ref{fig:RhombiMesh}.
This mesh is refined applying an anisotropic rule. In particular, at each step the mesh is refined by a factor $\alpha$ in the x-direction and by a factor $\alpha^2$ in the y-direction.
In Figure \ref{fig:convPoissonRhombi}, we present the convergence plots. We observe that the standard VEM method is not properly converging, while the proposed scheme exhibit the expected convergence behaviour.

\section{Conclusions}
We have presented a self-stabilized Virtual Element Method for the Poisson problem in mixed form. One of the main features of our approach is the employment of a projection operator over the gradients of harmonic polynomials of suitable degree. This choice alleviates the computational costs arising from the numerical quadrature. 
Despite the scheme is designed for arbitrary polygons, the theoretical analysis has been developed only for quadrilateral meshes.
The method convergence and stability have been computationally confirmed. Moreover, the numerical results show that our scheme is a valid alternative to the standard lowest-order mixed VEM.

A possible future development of the present study is the extension of the analysis to general polygonal meshes. 

\section*{Acknowledgements}
 The authors kindly acknowledge partial financial support by INdAM-GNCS projects 2022  CUP\_E55F2200027001. C.L. and M.V. kindly acknowledge partial financial support by PRIN 2017 (No. 201744KLJL) and PRIN 2020 (No. 20204LN5N5), funded by the Italian Ministry of Universities and Research (MUR).
 A.B and F.M. kindly acknowledge financial support provided by PNRR M4C2 project of CN00000013 National Centre for HPC, Big Data and Quantum Computing (HPC) CUP:E13C22000990001.
 A.B. kindly acknowledges partial financial support provided by INdAM-GNCS Projects 2023, MIUR project ``Dipartimenti di Eccellenza'' Programme (2018–2022)
    CUP:E11G18000350001 and by the PRIN 2020 project (No. 20204LN5N5\_003).

\bibliographystyle{plain}
 \bibliography{bibliografia-Mixed}
\end{document}